\documentclass[12pt,a4paper]{article}
\usepackage[latin1]{inputenc}
\usepackage[english]{babel}
\usepackage{amsmath}
\usepackage{amsfonts}
\usepackage{amssymb}
\usepackage{indentfirst}
\usepackage{dsfont}
\newtheorem{theorem}{Theorem}[section]
\newtheorem{lemma}[theorem]{Lemma}

\newtheorem{corollary}[theorem]{Corollary}

\title{A Threshold Regularization Method for Inverse Problems}
\date{}
\newcommand{\mykeywords}{Inverse problems; regularization; oracle inequalities; hard thresholding}
\newcommand{\mysubjclass}{
 62G05,
 62G08
 }
\author{Paul Rochet}
\begin{document}
\maketitle
\begin{flushleft}
{\noindent \normalfont\footnotesize Institut de Mathématiques de Toulouse, Université Paul Sabatier, 118 Route de Narbonne, 31062 Toulouse, France. Tel.: +33 561556790; fax: +33 561557599.\\ 
E-mail address: rochet@math.univ-toulouse.fr}
\end{flushleft}
\begin{abstract} A number of regularization methods for discrete inverse problems consist in considering weighted versions of the usual least square solution. However, these so-called filter methods are generally restricted to monotonic transformations, e.g. the Tikhonov regularization or the spectral cut-off. In this paper, we point out that in several cases, non-monotonic sequences of filters are more efficient. We study a regularization method that naturally extends the spectral cut-off procedure to non-monotonic sequences and provide several oracle inequalities, showing the method to be nearly optimal under mild assumptions. Then, we extend the method to inverse problems with noisy operator and provide efficiency results in a newly introduced conditional framework.
\end{abstract}
$ $ \vskip .1in
\noindent \textbf{Keywords}: \mykeywords\\
\noindent \textbf{Subject Class. MSC-2000} : \mysubjclass

\section{Introduction}

We are interested in recovering an unobservable signal $x_0$, based on noisy observations of the image of $ x_0$ through a linear operator $A$. The observation $y$ satisfies the following relation
$$ y(t) = A x_0 (t) + \varepsilon(t),$$
where $\varepsilon(.)$ is a random process representing the noise. This problem is studied in \cite{Cavalier00oracleinequalities}, \cite{MR1241596}, \cite{MR2440445} and in many applied fields such as medical imaging in \cite{MR1847845} or seismography in \cite{MR966773} for instance. When the measured signal is only available at a finite number of points $t_1,...,t_n$, the operator $A$ must be replaced by a discrete version $A_n: x \mapsto (A x(t_1),...,A x(t_n))'$, leading to a discrete linear model
$$ y = A_n x_0 + \varepsilon, $$
with $y \in \mathbb R^n$. Difficulties in estimating $x_0$ occur when the problem is \textit{ill-posed}, in the sense that small perturbations in the observations induce large changes in the solution. This is caused by an ill-conditioning of the operator $A_n$, reflected by a fast decay of its spectral values $b_i$. In such problems, the least square solution, although having a small bias, is generally inefficient due to a too large variance. Hence, \textit{regularization} of the problem is required to improve the estimation. A large number of regularization methods are based on considering weighted versions of the least square estimator. The idea is to allocate low weights $\lambda_i$, or \textit{filters}, to the least square coefficients that are highly contaminated with noise, thus reducing the variance, at the cost of increasing the bias at the same time. The most famous filter-based method is arguably the one due to Tikhonov (see \cite{MR0455365}), where a collection of filters is indirectly obtained via a minimization procedure with $\ell^2$ penalization. Tikhonov filters are entirely determined by a parameter $\tau$ that controls the balance between the minimization of the $\ell^2$ norm of the estimator and the residual.

Another well spread filter method that will be given a particular interest in this paper, is the \textit{spectral cut-off} discussed in \cite{MR2361904}, \cite{MR1408680} and \cite{MR916729}. One simply considers a truncated version of the least square solution, where all coefficients corresponding to arbitrarily small eigenvalues are removed. Thus, spectral cut-off is associated to binary filters $\lambda_i$, equal to $1$ if the corresponding eigenvalue $b_i$ exceeds in absolute value a certain threshold $\tau$, and $0$ otherwise. 

A common feature of spectral cut-off and Tikhonov regularization is the predetermined nature of the filters $\lambda_i$, defined in each case as a fixed non-decreasing function $f(\tau,.)$ of the eigenvalues $b_i^2$, and where only the parameter $\tau$ is allowed to depend on the observations. However, in many situations, non-monotonic sequences of filters may provide a more efficient estimation of $x_0$. Actually, optimal values for $\lambda_i$ generally depend on both the noise level, which is determined by the eigenvalue $b_i$, and the component, say $x_i$, of $x_0$ in the direction associated to $b_i$. A restriction to monotonic collections of filters turns out to be inefficient in situations where the coefficients $x_i $ are uncorrelated to the spectral values $b_i$ of the operator $A_n$.
 
Regularization methods involving more general classes of filters have also been treated in the literature. In \cite{Cavalier00oracleinequalities}, the authors study a general procedure known as \textit{unbiased risk estimation}, that applies to arbitrary classes of filters, dealing in particular with non-monotonic collections. However, their general framework concerning the class of estimators requires in return additional regularity assumptions which we intend to relax in this paper. We focus on a specific class of projection estimators that extends the spectral cut-off to non-monotonic collections of filters. Precisely, we consider the collection of unrestricted binary filters $\lambda_i \in \{0,1 \}$. The computation of the estimator relies on the choice of a proper set of coefficients $m \subset \{ 1,...,n \}$, which considerably increases the number of possibilities compared to the spectral cut-off procedure. We show this method to satisfy a non-asymptotic exact oracle inequality, when the oracle is computed in the class of binary filters. Moreover, we show our estimator to nearly achieve the rate of convergence of the best linear estimator in the maximal class of filters, i.e. when no restriction is made on $\lambda_i$.  

It many actual situations, the operator $A_n$ is not known precisely and only an approximation of it is available. Regularization of inverse problems with approximate operator is studied in \cite{MR2158113}, \cite{MR1872847} and \cite{MR2387973}. In this paper, we tackle the problem of estimating $x_0$ in the situation where we observe independently a noisy version $\hat b_i$ of each eigenvalue $b_i$. We consider a new framework where the observations $\hat b_i$ are made once and for all, and are seen as non-random. We provide a bound on the conditional risk of the estimator, given the values of $\hat b_i$, in the form of a conditional oracle inequality. 

The paper is organized as follows. We introduce the problem in Section \ref{thresholdsec2}. We define our estimator in Section \ref{thresholdsec3}, and provide two types of oracle inequalities. Section \ref{thresholdsec4} is devoted to an application of the method to inverse problems with noisy operators. The proofs of the results are postponed to the Appendix.

\section{Problem setting}\label{thresholdsec2}
Let $(\mathcal X, \Vert . \Vert)$ be a Hilbert space and $A_n: \mathcal X \to \mathbb R^n$ ($n >2$) a linear operator. We want to recover an unknown signal $x_0\in \mathcal X$ based on the indirect observations
\begin{equation}\label{modeldiscrete} y = A_n x_0 + \varepsilon,
\end{equation}
where $\varepsilon$ is a random noise vector. We assume that $\varepsilon$ is centered with covariance matrix $\sigma^2 I$, where $I$ denotes the identity matrix. We endow $\mathbb R^n$ with the scalar product $\left\langle u,v \right\rangle_n = n^{-1} \sum_{i=1}^n u_i v_i$ and the associated norm $\Vert . \Vert_n$ and we note $A_n^*: \mathbb R^n \to \mathcal X$ the adjoint of $A_n$. Let $\mathcal K_n$ be the kernel of $A_n$ and $\mathcal K_n^\perp$ its orthogonal in $\mathcal X$ which we assume to be of dimension $n$. The surjectivity of $A_n$ ensures that the observation $y$ provides information in all directions. If this condition is not met, one may simply reduce the dimension of the image in order to make $A_n$ surjective. 

The efficiency of the estimator relies first of all on the accuracy of the discrete operator $A_n$ and how "close" it is to the true value $A$. The convergence of the estimator towards $x_0$ is subject to the condition that the distance of $x_0$ to the set $\mathcal K_n^\perp$ tends to $0$, which is reflected by a proper asymptotic behavior of the design $t_1,...,t_n$. This aspect is not discussed here, we consider a framework where we have no control over the design $t_1,...,t_n$ and we focus on the convergence of the estimator towards the projection $x^\dagger$.

Let $\{b_i; \phi_i, \psi_i\}_{i=1,...,n}$ be a singular system for the linear operator $A_n$, that is, $A_n \phi_i = b_i \psi_i$ and $A_n^* \psi_i = b_i \phi_i$ and $b_1^2 \geq ... \geq b_n^2 > 0$ are the ordered non-zero eigenvalues of the self-adjoint operator $A_n^* A_n$. The $\phi_i$'s (resp. $\psi_i$'s) form an orthonormal system of $\mathcal K_n^\perp$ (resp. $\mathbb R^n$).\\

In this framework, the available information on $x_0$ consists in a noisy version of $A_n x_0$. As a result, estimating the part of $x_0$ lying in $\mathcal K_n$ is impossible, based only on the observations. The best approximation of $x_0$ one can get without prior information is the orthogonal projection of $x_0$ onto $\mathcal K_n^\perp$. This projection, noted $x^\dagger$, is called \textit{best approximate solution} and is obtained as the image of $A_n x_0$ through the generalized Moore-Penrose inverse operator $A_n^\dagger = (A_n^* A_n)^{\dagger} A_n^*$, where $(A_n^* A_n)^{\dagger}$ denotes the inverse of $A_n^* A_n$, restricted to $\mathcal K_n^\perp$. By construction, the generalized Moore-Penrose inverse $A_n^\dagger$ can also be defined as the operator for which $\{b_i^{-1}; \psi_i, \phi_i \}_{i=1,...,n}$ is a singular system. We refer to \cite{MR1408680} for further details.\\

Searching for a solution in the subspace $\mathcal K_n^\perp$ allows to reduce the number of regressors to $n$. Then, estimating $x^\dagger$ can be made using a classical linear regression framework where the number of regressors is equal to the dimension of the observation. Decomposing the observation in the singular basis $\{\psi_i\}_{i=1,...,n}$ leads to the following model
$$ y_i = b_i x_i + \varepsilon_i, i =1,...,n,  $$
where we set $y_i = \left\langle y, \psi_i \right\rangle_n$, $\varepsilon_i = \left\langle \varepsilon, \psi_i \right\rangle_n$ and $x_i = \left\langle x_0, \phi_i \right\rangle$. It now suffices to divide each term by the known singular value $b_i$ to observe the coefficient $x_i$, up to a noise term $\eta_i := b_i^{-1} \varepsilon_i$. Equivalently, this is obtained by applying the Moore-Penrose inverse $A_n^\dagger$ in the model \eqref{modeldiscrete}. We thus consider the function $y^\dagger = A_n^\dagger y \in \mathcal K_n^\perp$, defined as the inverse image of $y$ through $A_n$ with minimal norm. Identifying $y^\dagger$ with the vector of its coefficients $y_i^\dagger = b_i^{-1} y_i$ in the basis $\{\phi_i \}_{i=1,...,n}$, we obtain
\begin{equation}\label{modelr2} y^\dagger_i  = x_i + \eta_i, \ i=1,...,n. \end{equation}
The covariance matrix of the noise $\eta = (\eta_1,...,\eta_n)'$ is diagonal in this model, as we have $\mathbb E(\eta_i \eta_j) =  n^{-1} b_i^{-1} b_j^{-1} \sigma^2  \left\langle \psi_i,\psi_j\right\rangle_n$ which is null for all $i \neq j$ and equal to $\sigma_i^2 := \sigma^2 b_i^{-2} /n$ if $i =j$. Thus, the model can be interpreted as a linear regression model with heteroscedastic noises, the variances $\sigma_i^2$ being inversely proportional to the eigenvalues $b_i^2$. In the case where $\varepsilon$ in the original model \eqref{modeldiscrete} is Gaussian with distribution $\mathcal N(0,\sigma^2 I)$, the noises $\eta_i$ remain Gaussian in \eqref{modelr2}.\\

This representation points out the effect of the decay of the singular values $b_i$ on the noise level, making the problem ill-posed. To control the noise with a too large variance $\sigma_i^2$, a solution is to consider weighted versions of $y^\dagger$. For some filter $\lambda = (\lambda_1,...,\lambda_n)'$, note $\hat x(\lambda) \in \mathcal K_n^\perp$ the function defined by $\left\langle \hat x(\lambda), \phi_i \right\rangle = \lambda_i y_i^\dagger$ for $i=1,...,n$. Filter-based methods aim to cancel out the high frequency noises by allocating low weights to the components $y_i^\dagger$ corresponding to small singular values. A widely used example is the Tikhonov regularization, with weights of the form $\lambda_i = (1+ \tau \sigma_i^2)^{-1}$ for some $\tau >0$. The Tikhonov solution can be expressed as the minimizer of the functional
$$ \Vert y - A_n x \Vert^2 + \tau \Vert x \Vert^2, \ x \in \mathcal X,  $$
which makes the method particularly convenient in cases where the SVD of $A_n^* A_n$ or the coefficients $y_i^\dagger$ are not easily computable. We refer to \cite{MR2421941} and \cite{MR0455365} for further details.

Another common filter-based method is the \textit{truncated singular value decomposition} or \textit{spectral cut-off} studied in \cite{MR2361904}, \cite{MR1408680} and \cite{MR916729}. An estimator of $x_0$ is obtained as a truncated version of $y^\dagger$, where all coefficient $y_i^\dagger$ corresponding to arbitrarily small singular values are replaced by $0$. This approach can be viewed as a principal component analysis, where only the highly explanatory directions are selected. The spectral cut-off estimator is associated to filter factors of the form $\lambda_i = \mathds 1  \{i \leq k \}$, where $\mathds 1  \{ . \}$ denotes the indicator function and $k$ is a bandwidth to be determined. Data-driven methods for selecting suitable values of $k$ are discussed in \cite{MR2421941}, \cite{MR2283712}, \cite{MR916729}, \cite{MR0334486} and  \cite{MR516385}. \\

A natural way to generalize the spectral cut-off procedure is to enlarge the class of estimators by considering non-ordered truncated versions of $y^\dagger$, as made in \cite{MR2440445}, \cite{MR2426106} or \cite{MR2741964} (see also Examples 1 and 2 in \cite{Cavalier00oracleinequalities}). This approach reduces to a model selection issue where each model is identified with a set of indices $m \subset \{1,...,n \}$. Precisely, for $m$ a given model, define $\hat x_m \in \mathcal K_n^\perp$ as the orthogonal projection of $y^\dagger$ onto $\mathcal X_m := \text{span} \{ \phi_i, i \in m \}$, that is, $\hat x_m$ satisfies 
 $$\left\langle \hat x_m, \phi_i \right\rangle = \left\{ \begin{array}{c} y^\dagger_i \ \text{ if } \ i \in m, \\ 
 \ 0 \ \ \text{ otherwise.} \end{array} \right.$$ 
The objective is to find a model $m$ that makes the expected risk $ \mathbb E \Vert \hat x_m - x_0 \Vert^2$ small. The computation of the estimator no longer relies on the choice of one parameter $k \in \{1,...,n \}$ as for spectral cut-off, but on the choice of a set of indices $m \subset \{1,...,n \}$, which increases the number of possibilities. In particular, this approach allows non-monotonic collections of filters that may perform better than decreasing sequences obtained by spectral cut-off. To see this, write the bias-variance decomposition of the estimator $\hat x_m$ for a deterministic model $m$:
$$  \mathbb E \Vert \hat x_m - x_0 \Vert^2 = \mathbb E \Vert x_0 - x^\dagger \Vert^2 + \sum_{ i \notin m} x_i^2 +  \sum_{ i \in m} \sigma_i^2.  $$
In these settings, it appears that in order to minimize the risk, best is to select indices $i$ for which the component $x_i^2$ is larger than the noise level $\sigma_i^2$. A proper choice of filter should depend on both the variance $\sigma_i^2$ and the coefficient $x_i^2$. Consequently, the resulting sequence $\{ \lambda_i \}_{i=1,...,n}$ has no reason of being a decreasing function of $\sigma_i^2$ if some coefficients $x_i^2$ are large enough to compensate for a large variance.

\section{Non-ordered variable selection}\label{thresholdsec3}
\subsection{Threshold regularization}
The construction of the estimator by non-ordered variable selection reduces to finding a proper set $m$. Following the discrepancy principle, an optimal value for $m$ (minimizing the risk) is obtained by keeping small simultaneously the bias term $\sum_{ i \notin m} x_i^2$ and the variance term $\sum_{ i \in m} \sigma_i^2$ in the expression of the risk $\mathbb E \Vert \hat x_m - x_0 \Vert^2$. Following the previous argument, a minimizer of the risk $  \mathbb E \Vert \hat x_m - x_0 \Vert^2$ is obtained by selecting only the indices $i$ for which the coefficient $x_i^2$ is larger than the noise level $\sigma_i^2$. An optimal model is thus given by  $ m^* := \left\{ i: \ x_i^2 \geq \sigma_i^2 \right\}$. The coefficients $x_i$ being unknown to the practitioner, the optimal set $m^*$ can not be computed in practical cases. For this reason it will be referred to as an \textit{oracle}. \\

We shall now provide a model $\widehat m$ constructed from the available information, that mimics the oracle $m^*$. Fixing a threshold on the coefficients $x_i$ being impossible, we propose to use a threshold on the coefficients $y^\dagger_i$. Precisely, consider the set
$$\widehat m = \left\{ i:  y^{\dagger 2}_i \geq 4 \sigma_i^2 \mu_i \right\} ,$$ 
for $\{ \mu_i \}_{i=1,...,n}$ a sequence of positive parameters to be chosen. Obviously, the behavior of the resulting estimator $\hat x_{\widehat m}$ relies on the choice of the sequence $\{ \mu_i \}_{i=1,...,n}$: the larger the $\mu_i$'s, the more sparse is $\hat x_{\widehat m}$. It must be chosen so that the resulting set $\widehat m$ contains only the indices $i$ for which the noise level is small compared to the actual value of $x_i$. Although, the only knowledge of the observations $y_i^\dagger$ and the variances $\sigma_i^2$ makes it a difficult task. \\

There exist general filter-based methods that can be applied to arbitrary classes of filter estimators. One example is the \textit{unbiased risk estimation} discussed in \cite{Cavalier00oracleinequalities}, which defines an estimator of $x_0$ via the minimization of an unbiased estimation of the risk, over an arbitrary set of filters. When restricted to the class of binary filters $\lambda_i \in \{0,1 \}$, unbiased risk estimation reduces to minimizing over $\mathcal M$ the criterion
$$ m \mapsto \Vert y^\dagger - \hat x_m \Vert^2 + 2 \sum_{i \in m} \sigma_i^2.  $$
The minimum can be shown to be reached for the set $m = \{i : y_i^{\dagger 2} \geq 2 \sigma_i^2 \}$, which corresponds to taking $\mu_i = 1/2$ in our method. This choice is shown to be asymptotically efficient in Proposition 2 in \cite{Cavalier00oracleinequalities}, although additional restrictions are made on the $\lambda_i$'s which we intend to relax here in an asymptotic framework. If these conditions are not met, the accuracy of the choice $\mu_i = 1/2$ is not clear. We investigate in the next section a different choice for $\mu_i$ which turns out to be nearly optimal in a general framework.\\

\noindent In a general point of view, the estimator $\hat x_{\widehat m}$ can be obtained via a minimization procedure, using a BIC-type criterion for heteroscedastic models,
$$ \hat x_{\widehat m} = \underset{x \in \mathcal X}{\text{ argmin }} \left\{ \Vert y^\dagger - x \Vert^2 + 4  \sum_{i =1}^n \sigma_i^2 \mu_i \mathds 1  \{ \left\langle x, \phi_i \right\rangle \neq 0 \} \right\}. $$
In a certain way, this can be seen as a hard-thresholding version of the estimator considered in \cite{MR2741964}, obtained with a $\ell^1$ penalty. However, expressing the estimator as the solution to a minimization equation does not ease the computation. The method requires in any case calculation of the SVD of $A_n^* A_n$ and the coefficients $y_i^\dagger$, which may be computationally expensive. On the other hand, the computation of the estimator is simple once the decomposition of $y^\dagger$ in the SVD of $A_n^* A_n$ is known, as it suffices to compare each coefficient $y_i^\dagger$ to the threshold $4 \sigma_i^2 \mu_i$.

\subsection{Oracle inequalities}\label{threshold31}
In the definition of $\widehat m$, the choice of the parameters $\mu_i$ is crucial. Too large values of $\mu_i$ will result in an under-adjustment, keeping too few relevant components $y_i^\dagger$ to estimate $x_0$. On the contrary, a small value of $\mu_i$ increases the probability of selecting a component $y_i^\dagger$ that is highly affected with noise. Thus, it is essential to find a good balance between these two types of errors. In the next theorem, we provide a nearly optimal choice for the parameters $\mu_i$, under the condition that $\varepsilon$ has finite exponential moments.\\ 

\noindent For $i=1,...,n$, note $\gamma_i: = \eta_i^2/\sigma_i^2 = n \varepsilon_i^2/\sigma^2$. We make the following assumption.\\

 A1. There exist $K, \beta >0$ such that $ \forall t>0, \forall i=1,...,n,  \ \mathbb P(\gamma_i > t) \leq K e^{-t/\beta}$.\\
 
\noindent In a Gaussian model, the $\gamma_i$'s have $\chi^2$ distribution with one degree of freedom. The condition A1 holds for any $\beta >2$, taking $K = \sqrt{1-2/\beta}$.

\begin{theorem}\label{oraclekbeta} Assume that the condition \emph{A1} holds. Set $ \mu_i = \emph{max}\{\beta \log(n^2 \sigma_i^2),0 \}$, the estimator $\hat x_{\widehat m}$ satisfies
$$ \mathbb E \Vert \hat x_{\widehat m} - x^\dagger \Vert^2 \leq \mathbb E \Vert \hat x_{m^*} - x^\dagger \Vert^2  + ( K_1 \log n + K_2) \sum_{i \in m^*} \sigma_i^2 + \frac {K_3} n ,   $$
with $K_1 =12 \beta $, $K_2 = 2 + \beta \log \Vert x^\dagger \Vert^2$ and $K_3 = 2 K \beta$.
\end{theorem}

Remark 1. This theorem establishes a non-asymptotic oracle inequality with exact constant. The residual term is similar to that in Corollary 1 in \cite{MR2440445}. The fact that the term $\Vert x^\dagger \Vert$ depends on $n$ is not problematic here as it can in any case be bounded by the norm of $x_0$. \\

Remark 2. The method requires knowledge of the operator $A_n$, the variance $\sigma^2$ and the constant $\beta$ in the condition A1. Note however that knowing the constant $K$ is not necessary to build the estimator.  \\

Remark 3. The set $\widehat m$ contains all indices $i$ for which $\sigma_i^2 \leq 1/n^2$, as we have in this case $\mu_i =0$. This suggests that the error caused by wrongfully selecting indices $i$ for which the variance is smaller than $1/n^2$ is negligible, regardless of the value of $y_i^\dagger$.\\

Remark 4. In an asymptotic concern, the accuracy of the result stated in Theorem \ref{oraclekbeta} relies on the convergence rate of the residual term to zero, compared to the risk of the oracle. The residual term $\sum_{i \in m^*} \sigma_i^{2} $ is actually the variance term in the bias-variance decomposition of $\hat x_{m^*}$, and therefore, it is bounded by the risk of the oracle. As a result, the estimator $\hat x_{\widehat m}$ is shown to reach at least the same rate of convergence as the oracle, up to a logarithmic term, which warrants good adaptivity properties. The logarithmic term vanishes in the convergence rate if the bias term $\sum_{i \notin m^*} x_i^{2} $ dominates in the risk of the $\hat x_{m^*}$. Precisely, the oracle inequality is asymptotically exact as soon as the residual term $\log n \sum_{i \in m^*} \sigma_i^{2}$ is negligible compared to the bias term $\sum_{i \notin m^*} x_i^{2}$. In this case, it follows from Theorem \ref{oraclekbeta} that
$$ \mathbb E \Vert \hat x_{\widehat m} - x^\dagger \Vert^2 = (1 + o(1)) \ \mathbb E \Vert \hat x_{m^*} - x^\dagger \Vert^2.   $$
Of course, this condition is hard to verify in practice and assuming it is true reduces to make strong regularity assumptions on the asymptotic behavior of $x_0$ and $A_n$. In a non-asymptotic framework, the theorem warrants that the estimator $\hat x_{\widehat m}$ is close to the oracle as soon as the variance term $ \sum_{i \in m^*} \sigma_i^{2}$ is small compared to the bias term $\sum_{i \notin m^*} x_i^{2}$ in the bias-variance decomposition of the oracle. \\
 
The estimator $\hat x_{\widehat m}$ being built using binary filters $\lambda_i \in \{0,1 \}$, it is natural to measure its efficiency by comparing its risk to that of the best linear estimator in this class. Nevertheless, we see in the next corollary that a similar oracle inequality holds if we consider the oracle in the maximal class of filters, that is, allowing the $\lambda_i$'s to take any real value.

\begin{corollary} Assume that the condition \emph{A1} holds, the estimator $\hat x_{\widehat m}$ of Theorem \ref{oraclekbeta} satisfies
$$ \mathbb E \Vert \hat x_{\widehat m} - x^\dagger \Vert^2 \leq K_4 \log n \ \inf_{\lambda \in \mathbb R^n} \mathbb E \Vert \hat x(\lambda) - x^\dagger \Vert^2 + \frac{K_5} n, $$
for some constants $K_4, K_5$ independent of $n$.
\end{corollary}

\noindent This result is a straightforward consequence of Lemma \ref{majoracles} in the Appendix, where it is shown that the oracle in the class of binary filters $\lambda_i \in \{0,1 \}$ achieves the same rate of convergence up to a factor $2$, as the best filter estimator obtained with non-random values of $\lambda$. The class of unrestricted binary filters leads to a simple solution while it induces a slight loss of efficiency compared to the maximal class. \\

Interest of oracles lies in the fact that the best estimator in a given class will often reach the optimal rate of convergence. In many situations, comparing the risk of the estimator to that of an oracle might be sufficient to deduce optimality results, as well as adaptivity properties, as discussed in \cite{MR2421941}. In the literature of inverse problems, rates of convergence of oracles are obtained under regularity conditions on the map $x_0$ and the spectrum of $A_n$. These conditions can be gathered into a single assumption, generally referred to as \textit{source condition}, relating the behavior of $x_0$ to the regularity of the operator $A_n$ (see for instance \cite{MR2361904}, \cite{MR1408680} or \cite{MR2434334}). Another point of view widely adopted in the literature is the \textit{minimax} approach (see \cite{MR2421941}), aiming to determine the behavior the worst possible value of $x_0$ in a given class of functions. Typically, the condition can be a polynomial decay of the coefficients $x_i$, which reduces to assuming that $x_0$ lies in the unit ball in a proper Besov space. For rates of convergence with a minimax approach, we refer to \cite{MR5193414}, \cite{MR1394051} and \cite{MR2027536}. In our framework, rates of convergence for $\hat x_{\widehat m}$ can be deduced from Theorem 2 in \cite{MR2440445}, under a polynomial decay of the coefficients $x_i$ and the eigenvalues $b_i$.

\section{Regularization with unknown operator}\label{thresholdsec4}
In many actual situations, the operator $A_n$ is not precisely known. In this section, we consider the framework where the operator $A_n$ is observed independently from $y$. This situation is treated in \cite{MR2158113}, \cite{MR1872847} or \cite{MR2387973}. The method discussed in the previous section does not apply for such problems since it requires complete knowledge of the operator $A_n$.\\ 
As in \cite{MR2158113}, we assume that the eigenvectors $\phi_i$ and $\psi_i$ are known. This seemingly strong assumption is actually met in many situations, for instance if the problem involves convolution or differential operators which can be decomposed in Fourier basis (see also the examples in \cite{MR2421941}). Thus, only the eigenvalues $b_i$ are unknown and we assume they are observed independently of $y$, with a centered noise $\xi_i$ with known variance $s^2>0$:
$$  \hat b_i = b_i + \xi_i, \ i =1,...,n.   $$
The method discussed in this paper is different according to whether the eigenvalues are known exactly or observed with a noise. Thus, we need to assume here that $s$ is positive and the known operator framework can not be seen as a particular case. Moreover, we assume the $\xi_i$'s are independent and satisfy the two following conditions.\\

 A2. There exist $K', \beta' >0$ such that $ \forall t>0, \forall i=1,...,n, \ \mathbb P( \xi_i^2/s^2 > t) \leq K' e^{-t/\beta'}$.\\

 A3. There exist $C, \alpha >0$ such that $  \forall i=1,...,n, \ \min \{ \mathbb P( \xi_i < - \alpha s),  \mathbb P( \xi_i > \alpha s) \} \geq C$.\\

\noindent As discussed previously, the condition A2 means that that the $\xi_i$'s have finite exponential moments. The condition A3 is hardly restrictive, and is fulfilled for instance as soon as the $\xi_i$'s are identically distributed. As we shall see in the sequel, the method requires knowledge of the constant $\alpha$ (or at least an upper bound for it), but no information on the constants $\beta'$, $K'$ or $C$ is needed to build the estimator. \\

\noindent Knowing the eigenvectors of $A_n^* A_n$ allows us to write the model in the form
$$ y_i = b_i x_i + \varepsilon_i, i =1,...,n.  $$
In our framework where the actual eigenvalues $b_i$ are unknown, a natural estimator of each component $x_i$ is obtained by $\tilde y_i = \hat b_i^{-1} y_i$, provided that $\hat b_i \neq 0$. However, it is clear that this estimate is not satisfactory if $\hat b_i$ is far from the true value (consider for instance the extreme case where $\hat b_i=0$ or if $\hat b_i$ and $b_i$ are of opposite signs). Actually, the naive estimator $\hat b_i^{-1}$ can not be used efficiently to estimate $b_i^{-1}$ because it may have an infinite variance. In \cite{MR2158113}, the authors fix a threshold $w$ the estimate can not exceed and consider an estimator of $b_i^{-1}$ equal to $\hat b_i^{-1}$ if $\vert \hat b_i \vert  > 1/w$ and null otherwise. As we will see below, we use the same idea here, although the threshold fixed on the $\hat b_i$'s is implicitly part of the variable selection process.  

We can reasonably assume that null values of $\hat b_i$ do not provide any relevant information and can not be used to estimate $x_0$. Thus, to avoid considering trivial situations, we assume that all $\hat b_i$ are non-zero. In all generality, the $\tilde y_i$'s can be viewed as noisy observations of $x_i$ by writing
$$ \tilde y_i = x_i + \tilde \eta_i, \ i=1,...,n,    $$
with $\tilde y_i = \hat b_i^{-1} \left\langle y,\psi_i \right\rangle_n$ and $\tilde \eta_i = \hat b_i^{-1} (\varepsilon_i - \xi_i x_i)$, where we recall $\varepsilon_i = \left\langle  \varepsilon, \psi_i \right\rangle_n$. As in the previous section, we propose a threshold procedure to filter out the observations $\tilde y_i$ that are potentially highly contaminated with noise. Here, the noise $\tilde \eta_i$ is more difficult to deal with because it depends on the unknown coefficient $x_i$. 

Our objective is to find an optimal variable selection criterion conditionally to the $\hat b_i$'s. In order to do so, we consider a framework where the $\hat b_i$'s are observed once and for all, and are treated as non-random. Thus, we define as an oracle, a model $m^*_\xi$ minimizing the conditional risk $ \mathbb E_\xi \Vert \hat x_m - x^\dagger \Vert^2$, where $\mathbb E_\xi(.)$ denotes the expectation knowing $\xi=(\xi_1,...,\xi_n)'$. Following a similar argument as in the previous section, a model minimizing the conditional risk contains only the indices $i$ for which the coefficient $x_i^2$ is larger than the noise level. Hence, we may define $m^*_\xi = \{ i: x_i^2 > \mathbb E_\xi (\tilde \eta_i^2)\}$. A notable difference here is that the noise $\tilde \eta_i$ actually depends on the value $x_i$. Let $\hat \sigma_i^2 = n^{-1} \hat b_i^{-2} \sigma^2$, we can calculate the conditional expectation of $\tilde \eta_i^2$, given by $\mathbb E_\xi (\tilde \eta_i^2) = \hat \sigma_i^2 + \hat b_i^{-2} \xi_i^2 x_i^2$. After simplifications, it appears that the optimal model conditionally to the $\xi_i$'s can be expressed in the two following explicit forms
$$ m^*_\xi = \left\{ i: 2 \vert \hat b_i \vert > \frac{\sigma^2}{n \vert b_i \vert x_i^2} + \vert b_i\vert  \right\} = \left\{ i:  x_i^2 > \frac{\sigma^2}{n(\hat b_i^2 - \xi_i^2)}, \ \vert \hat b_i \vert > \frac{ \vert b_i \vert} 2  \right\}.$$
In the first expression, we see that the oracle selects indices $i$ for which the observation $\hat b_i$ exceeds a certain value depending on both $x_i$ and $b_i$. Interestingly, components $\tilde y_i$ corresponding to observations $\hat b_i$ smaller than half the true eigenvalue $b_i$ are not selected in the oracle, regardless of the coefficient $x_i$. Here again, the optimal model $ m^*_\xi$ can not be used in practical cases since it involves the unknown values $x_i$ and $\xi_i$. We can only try to mimic the optimal threshold, based on the observations $\tilde y_i$ and $\hat b_i$. Consider the set
$$ \widehat m_\xi = \left\{ i: \tilde y_i^2  > 8\hat \sigma_i^2 \nu_i  , \ \vert \hat b_i \vert >  \alpha s \right\},$$ 
where $\{ \nu_i \}_{i=1,...,n}$ are parameters to be chosen and $\alpha$ is the constant defined in A3. With this definition, only the indices for which the observation $\hat b_i$ is larger than a certain value, namely $\alpha s$, are selected. This conveys the idea discussed in \cite{MR2158113}, that when $b_i$ is small compared to the noise level, the observation $\hat b_i$ is potentially mainly noise. Remark however that in \cite{MR2158113}, the lower limit for the observed eigenvalues is $s \log^2 (1/s)$, while in our method, it is chosen of the same order as the standard deviation $s$. \\

\noindent Define the set $M =  \{ i:  \vert b_i \vert < 2\alpha s  \}$.

\begin{theorem}\label{oraclerandomth} Assume that the condition \emph{A1} holds. The threshold estimator obtained with $\nu_i = \emph{max}\{ \beta \log(n^2 \hat \sigma_i^2),0 \} $ satisfies,
$$ \mathbb E_\xi \Vert \hat x_{\widehat m_\xi} - x^\dagger \Vert^2 \leq \left(K_1' \log n + K'_2 \right) \mathbb E_\xi \Vert \hat x_{m^*_\xi} - x^\dagger \Vert^2 + \sum_{i \in M} x_i^2 + \kappa(\xi), $$  
with $K_1' = \emph{max} \{ 18 \beta, 4 \alpha^{-2} \beta' \} $, $K'_2 = \emph{max} \{ 9(\beta \log \Vert x^\dagger \Vert^2 +1), 1 \}$, and 
$$  \kappa(\xi)= \frac{4 K \beta} n + 4 \sum_{i \notin m^*_\xi} \frac{\xi_i^2 x_i^2}{\alpha^{2} s^{2}}  \mathds 1  \{\xi_i^2 >  s^2 \beta' \log n  \}.$$
Moreover, if \emph{A2} holds, $\mathbb E(\kappa(\xi)) = O(n^{-1} \log n).$
\end{theorem}
The main interest of this result lies in the fact that it provides an oracle inequality, conditionally to the $\hat b_i$'s. In particular, the conditional oracle $ \hat x_{m^*_\xi}$ is more efficient than the estimator obtained by minimizing the expected risk $m \mapsto \mathbb E \Vert \hat x_m - x^\dagger \Vert^2$, since the optimal set $m^*_\xi$ is allowed to depend on the $\xi_i$'s. We see that the estimator $\hat x_{\widehat m_\xi}$ performs almost as well as the conditional oracle. Indeed, the residual term $  \kappa(\xi)$ is independent from $\xi$ with high probability, and its expectation is negligible under A2 as pointed out in the theorem. The non-random term $\sum_{i \in M} x_i^2$ is small if the eigenvalues $b_i$ are observed with a good precision, i.e. if the variance $s^2$ is small. Moreover, this term can be shown to be of the same order as the risk under the condition A3.

\begin{corollary}\label{randomopcor} If the conditions \emph{A1}, \emph{A2} and \emph{A3} hold, the threshold estimator defined in Theorem \ref{oraclerandomth} satisfies
 $$ \mathbb E \Vert \hat x_{\widehat m_\xi} - x^\dagger \Vert^2 \leq K'_4 \log n \ \mathbb E \Vert \hat x_{m^*_\xi} - x^\dagger \Vert^2 + \frac{K'_5 \log n} n, $$
for some constants $K'_4$ and $K'_5$ independent from $n$ and $s^2$.
\end{corollary}
With a noisy operator, we manage to provide an estimator that achieves the rate of convergence of the conditional oracle, regardless of the precision of the approximation of the spectrum of $A_n$. Indeed, the constants $K'_4$ and $K'_5$ in Corollary \ref{randomopcor} do not involve the variance $s^2$ of $\xi$. Actually, the variance only plays a role in the accuracy of the oracle. The result is non-asymptotic and requires no assumption on $s^2$.

\section{Appendix}
\subsection{Technical lemmas}
\begin{lemma}\label{21}  Assume the condition \emph{A1} holds. We have  
\begin{itemize} 
 \item $ \mathbb E \left( (\eta_i^2 - x_i^2) \mathds 1  {\{i \in \widehat m\}} \right) \leq 2 K \beta \sigma_i^2 e^{-\mu_i/\beta}.$\vspace{-0.2cm}
 \item $ \mathbb E \left((x_i^2 - \eta_i^2) \mathds 1  { \{i \notin \widehat m \}} \right) \leq \sigma_i^2 ( 6 \mu_i + 2)$.
\end{itemize}
\end{lemma}
\textit{Proof.} Using the inequality $(a+b)^2 \leq 2a^2 + 2 b^2$, we find that $ \eta_i^2 - x_i^2 \leq 2 \eta_i^2 - y_i^{\dagger 2}/2$. By definition of $\widehat m$, we get
$$ (\eta_i^2 - x_i^2) \mathds 1  {\{i \in \widehat m\}} \leq 2 \sigma_i^2 (\gamma_i - \mu_i)\mathds 1  {\{i \in \widehat m\}} \leq 2 \sigma_i^2 (\gamma_i - \mu_i)\mathds 1  {\{\gamma_i \geq \mu_i\}},   $$
where we used that $X \leq X \mathds 1  {\{ X \geq 0 \}}$. We finally obtain for all $i \notin m^*$,
$$ \mathbb E \left( (\eta_i^2 - x_i^2) \mathds 1  {\{i \in \widehat m\}} \right) \leq 2 \sigma_i^2 \int_0^{\infty} \mathbb P(\gamma_i \geq t + \mu_i) \ dt \leq 2 K \beta \sigma_i^2 e^{-\mu_i/\beta}, $$
as a consequence of A1. For the second part of the lemma, write $x_i^2-\eta_i^2 = y_i^{\dagger 2} - 2 \eta_i y_i^\dagger$ which is bounded by $ 3 y_i^{\dagger 2}/2 + 2 \eta_i^2$, using the inequality $2ab \leq 2a^2 + b^2/2$. This leads to
\begin{eqnarray*} \mathbb E \left((x_i^2-\eta_i^2) \mathds 1  { \{i \notin \widehat m \}}\right)  \leq  \sigma_i^2 (6 \mu_i  + 2 ).
\end{eqnarray*}

\begin{lemma}\label{majoracles} 
$$ \inf_{m \in \mathcal M} \mathbb E \Vert \hat x_m - x^\dagger \Vert^2 \leq 2  \inf_{\lambda \in \mathbb R^n} \mathbb E \Vert \hat x(\lambda) - x^\dagger \Vert^2.   $$
\end{lemma}

\noindent \textit{Proof.} The minimal values of the expected risks can be calculated explicitly in the two classes considered here. Minimizing over $\mathbb R^n$ the function $\lambda \mapsto \mathbb E \Vert \hat x(\lambda) - x^\dagger \Vert^2$, we find that the optimal value of $\lambda_i$ is reached for $ \lambda_i^* = x_i^2/( x_i^2 + \sigma_i^2)$. On the other hand, we know that $m \mapsto \mathbb E \Vert \hat x_m - x^\dagger \Vert^2$ reaches its minimum at $ m^* = \{ i: x_i^2 \geq \sigma^2_i \}$, yielding
$$ \inf_{\lambda \in \mathbb R^n} \mathbb E \Vert \hat x(\lambda) - x^\dagger \Vert^2 = \sum_{i=1}^n \frac{x_i^2 \sigma_i^2 }{ x_i^2 + \sigma_i^2} \ \ \text{ and } \ \inf_{m \in \mathcal M} \mathbb E \Vert \hat x_m - x^\dagger \Vert^2 = \sum_{i \in m^*} \sigma_i^2  + \sum_{i \in m^*} x_i^2 .  $$
By definition, if $i \in m^*$, $2x_i^2/( x_i^2 + \sigma_i^2) \geq 1$. In the same way, $2\sigma_i^2/( x_i^2 + \sigma_i^2) \geq 1$, for all $i \notin m^*$. We conclude by summing all the terms.\\

\begin{lemma}\label{lemmaunknownop} Assume the condition \emph{A1} holds. We have, for all $i=1,...,n$,  
\begin{itemize} 
\item $ \mathbb E_\xi \left( (\tilde \eta_i^2 - x_i^2) \mathds 1  {\{i \in \widehat m_\xi\}} \right) \leq 4 K \beta \ \hat \sigma_i^2 e^{- \nu_i/\beta}  + \dfrac{4 \xi_i^2 x_i^2}{\alpha^{2} s^{2}}$.
\item $  \mathbb E_\xi \left( (x_i^2- \tilde \eta_i^2) \mathds 1  { \{i \notin \widehat m_\xi \}} \right) \leq    9 \hat \sigma_i^2 \nu_i + 8 \mathbb E_\xi (\tilde \eta_i^2)  + x_i^2 \mathds 1  \{ \vert \hat b_i \vert \leq \alpha s \}$.
\end{itemize}
\end{lemma}

\noindent \textit{Proof.} Remark that $ \tilde \eta_i^2 = \hat b_i^{-2}(\varepsilon_i - \xi_i x_i)^2 \leq 2 \hat b_i^{-2} \varepsilon_i^2 + 2 \hat b_i^{-2} \xi_i^2 x_i^2$. Using that $x_i^2 \geq \tilde y_i^2/2 - \tilde \eta_i^2$, we deduce
$$ \tilde \eta_i^2 - x_i^2 \leq 4\hat b_i^{-2}\varepsilon_i^2 + 4 \hat b_i^{-2} \xi_i^2 x_i^2 - \frac{\tilde y_i^2} 2.$$ 
Writing $\widehat m_\xi = \{ \tilde y_i^2 > 8 \hat \sigma_i^2 \nu_i \} \cap \{ \vert \hat b_i \vert > \alpha s \} $, we find
$$ (\tilde \eta_i^2 - x_i^2) \mathds 1  {\{i \in \widehat m_\xi\}} \leq 4 \hat \sigma_i^2 (\gamma_i -  \nu_i) \mathds 1  {\{ \gamma_i \geq \nu_i \}} + 4\hat b_i^{-2} \xi_i^2 x_i^2 \mathds 1  \{ \vert \hat b_i \vert > \alpha s \},   $$
where we recall that $\gamma_i = n \varepsilon_i^2/\sigma^2$. Clearly, $\hat b_i^{-2} \mathds 1  \{ \vert \hat b_i \vert > \alpha s \} < \alpha^{-2} s^{-2}$ and the result follows using the condition A1. For the second part of the lemma, remark that the complement of $\widehat m_\xi$ is $\{ \tilde y_i^2 \leq 8 \hat \sigma_i^2 \nu_i, \ \vert \hat b_i \vert >  \alpha s \} \cup \{ \vert \hat b_i \vert \leq  \alpha s \}$. Using the inequality $x_i^2- \tilde \eta_i^2 \leq (1+\theta^{-1}) \tilde y_i^2 + \theta \tilde \eta_i^2$ for $\theta = 8$, we get
\begin{eqnarray*} (x_i^2- \tilde \eta_i^2) \mathds 1  { \{i \notin \widehat m_\xi \}} \leq 9 \hat \sigma_i^2 \nu_i + 8 \tilde \eta_i^2  + x_i^2 \mathds 1  \{ \vert \hat b_i \vert \leq \alpha s \}. 
\end{eqnarray*}

\begin{lemma}\label{lemmacondit} If \emph{A2} holds, we have
$$  \xi_i^2 \leq s^2 \beta' \log n  + \xi_i^2 \mathds 1  \{\xi_i^2 >  s^2 \beta' \log n  \}, $$
with $ \mathbb E  \left(\xi_i^2 \mathds 1  \{\xi_i^2 >  s^2 \beta' \log n  \} \right) = O (n^{-1}\log n )$.
\end{lemma}

\noindent \textit{Proof.} Write $ \xi_i^2 \leq s^2 \beta' \log n \ \mathds 1  \{\xi_i^2 \leq s^2 \beta' \log n \} + \xi_i^2 \mathds 1  \{\xi_i^2 >  s^2 \beta' \log n  \}$. To bound the first term, we use the crude inequality $\mathds 1  \{\xi_i^2 \leq  s^2  \beta' \log n \} \leq 1$. For the second term, we have as a consequence of A2,
\begin{eqnarray*} \mathbb E \left( \xi_i^2 \mathds 1  \{\xi_i^2 > s^2 \beta' \log n \} \right) & = & \int_0^\infty \mathbb P \left( \xi_i^2 \mathds 1  \{\xi_i^2/s^2 >  \beta' \log n \} > t \right) \ dt \\
& = & s^2 \beta' \log n \ \mathbb P( \xi_i^2/s^2 >  \beta' \log n) + s^2 \int_{\beta' \log n}^\infty \mathbb P( \xi_i^2/s^2  > t) \ dt\\
&  \leq & \dfrac{ K'\beta' s^2 (1+\log n )}{ n}.
\end{eqnarray*}

\subsection{Proofs}
\noindent \textbf{Proof of Theorem \ref{oraclekbeta}.} Write 
$$ \Vert \hat x_{\widehat m} - x_0 \Vert^2 = \Vert \hat x_{m^*} - x_0 \Vert^2 + \sum_{ i \notin m^*} (\eta_i^2 - x_i^2) \mathds 1  {\{i \in \widehat m\} } + \sum_{ i \in m^*} ( x_i^2 - \eta_i^2) \mathds 1  { \{ i \notin \widehat m \} }. $$
The objective is to bound the terms $ \mathbb E((\eta_i^2 - x_i^2) \mathds 1  {\{i \in \widehat m\} })$ and $\mathbb E(( x_i^2 - \eta_i^2) \mathds 1  { \{ i \notin \widehat m \}})$. First, assume that $\sigma_i^2 > 1/n^2$, i.e. $\mu_i =\beta \log \left(n^2 \sigma_i^2 \right)$. By Lemma \ref{21}, we know that
\begin{eqnarray*} \mathbb E \left((\eta_i^2 - x_i^2 ) \mathds 1  { \{i \in \widehat m \}} \right) \leq 2 K \beta \sigma_i^2 e^{-\mu_i/\beta} \leq \frac{2K \beta}{n^2}.
\end{eqnarray*}
The same bound holds if $\sigma_i^2 \leq 1/n^2$ with $\mu_i =0$, as a straight-forward consequence of Lemma \ref{21}. On the other hand, note that if $i \notin \widehat m$, then $\mu_i = \beta \log(n^2 \sigma_i^2)$. Lemma \ref{21} warrants
$$ \mathbb E \left( (x_i^2 - \eta_i^2) \mathds 1  { \{i \notin \widehat m \}} \right) \leq \sigma_i^2 \left(6 \beta \log(n^2 \sigma_i^2) + 2\right).  $$
Since $i \in m^*$, $\log(n^2 \sigma_i^2) \leq 2 \log n + \log \Vert x^\dagger \Vert^2$. We conclude by summing all the terms.\\

\noindent \textbf{Proof of Theorem \ref{oraclerandomth}.} The proof starts as in Theorem \ref{oraclekbeta}. We have
$$ \Vert \hat x_{\widehat m_\xi} - x^\dagger \Vert^2 =  \Vert \hat x_{m_\xi^*} - x^\dagger \Vert^2 + \sum_{ i \notin m_\xi^*} (\tilde \eta_i^2 - x_i^2) \mathds 1  {\{i \in \widehat m_\xi \} } + \sum_{ i \in m_\xi^*} ( x_i^2 - \tilde \eta_i^2) \mathds 1  { \{ i \notin \widehat m_\xi \} }, $$
and the objective is to bound the conditional expectation of each term separately. Using successively Lemma \ref{lemmaunknownop} and Lemma \ref{lemmacondit}, we get
$$ \mathbb E_\xi \left( (\tilde \eta_i^2 - x_i^2) \mathds 1  {\{i \in \widehat m_\xi\}} \right) \leq \frac{4 K \beta}{n^2}  + 4 \alpha^{-2} s^{-2} \xi_i^2 x_i^2 \leq \frac{ 4\beta' \log n}{\alpha^2} \ x_i^2 + \kappa_i(\xi),  $$
with 
$$ \kappa_i(\xi) = \frac{4 K \beta}{n^2} + \frac{4\xi_i^2 x_i^2}{\alpha^2 s^2} \mathds 1  \{\xi_i^2 >  s^2 \beta' \log n  \}.  $$
By Lemma \ref{lemmacondit}, we know that $\kappa(\xi) = \sum_{i \notin m^*_\xi} \kappa_i(\xi)$ is such that 
$$\mathbb E (\kappa(\xi)) \leq \frac{ 4 ( K \beta + 2  \alpha^{-2} K' \beta' \Vert x^\dagger \Vert^2 \log n  )}{n} = O \left(\frac{\log n} n \right).$$ 
On the other hand, Lemma \ref{lemmaunknownop} gives, for $\theta = 8$,
$$ \mathbb E_\xi \left( (x_i^2- \tilde \eta_i^2) \mathds 1  { \{i \notin \widehat m_\xi \}} \right) \leq   9  \hat \sigma_i^2  \nu_i + 8 \mathbb E_\xi (\tilde \eta_i^2)  + x_i^2 \mathds 1  \{ \vert \hat b_i \vert \leq \alpha s \}. $$
For all $i \in m^*_\xi$, we know that $\vert \hat b_i \vert \geq \vert b_i \vert /2$. Thus, if $i \in m^*_\xi$, $\mathds 1  \{ \vert \hat b_i \vert \leq \alpha s \} \leq \mathds 1  \{ i \in M \}$, where we recall $M = \{ i: \vert b_i \vert < 2 \alpha s \}$. We know also that, if $i \in m^*_\xi$, then $\hat \sigma_i^2 \leq x_i^2$. Thus, $\nu_i = \beta \log(n^2 \hat \sigma^2_i) \leq 2 \beta \log n + \beta \log \Vert x^\dagger \Vert^2$. Noticing that $\hat \sigma_i^2 \leq \mathbb E_\xi( \tilde \eta_i^2)$, we find
$$ \mathbb E_\xi \left( (x_i^2- \tilde \eta_i^2) \mathds 1  { \{i \notin \widehat m_\xi \}} \right) \leq (18 \beta \log n + 9\beta \log \Vert x^\dagger \Vert^2 + 8) \mathbb E_\xi( \tilde \eta_i^2) + x_i^2 \mathds 1  \{ i \in M \}.  $$
The result follows by summing all the term, using that the risk of the oracle $\hat x_{\widehat m_\xi}$ is
$$ \mathbb E_\xi \Vert \hat x_{m_\xi^*} - x^\dagger \Vert^2 = \sum_{ i \notin m_\xi^*} x_i^2  +   \sum_{ i \in m_\xi^*} \mathbb E_\xi( \tilde \eta_i^2). $$

\noindent \textbf{Proof of Corollary \ref{randomopcor}.} It suffices to show that the term $\sum_{i \in M}  x_i^{2 }$ is of the same order as the risk of the oracle. Write
\begin{eqnarray*} \mathbb E \Vert \hat x_{m_\xi^*} - x^\dagger \Vert^2 \geq \sum_{i=1}^n  x_i^2 \mathbb P (i \notin m^*_\xi)  \geq \sum_{i=1}^n x_i^2 \mathbb P ( \vert \hat b_i \vert \leq \vert b_i\vert /2 ).  
\end{eqnarray*}
For all $i \in M$, the probability $\mathbb P ( \vert \hat b_i \vert \leq \vert b_i\vert /2 )$ is greater than $C$ as a consequence of A3. We deduce $ \sum_{i \in M}  x_i^{2} \leq C^{-1} \mathbb E \Vert \hat x_{m_\xi^*} - x^\dagger \Vert^2$.

\bibliographystyle{plain}
\bibliography{threshold}

\end{document}